\documentclass[11pt]{amsart}

\usepackage{amssymb,amscd,amsthm,amsxtra}
\usepackage{latexsym}

\newtheorem{thm}{Theorem}[section]
\newtheorem{cor}[thm]{Corollary}

\newtheorem{prop}[thm]{Proposition}
\theoremstyle{definition}
\newcommand{\comment}[1]{}

\newtheorem{defn}[thm]{Definition}
\theoremstyle{remark}
\newtheorem{rem}[thm]{Remark}
\numberwithin{equation}{section}
\newcommand{\R}{\mathbb R}
\newcommand{\C}{\mathbb C}
\newcommand{\CN}{\mathcal N}
\begin{document}

\title[Improved interaction Morawetz inequalities for cubic
NLS on $\R^2$]{Improved interaction Morawetz inequalities for the cubic
  nonlinear Schr\"odinger equation on $\R^2$}
\author{J. Colliander}
\thanks{J.C. was supported in part by NSERC grant RGP250233-03.}
\address{Department of Mathematics, University of Toronto, Toronto, ON, Canada M5S 2E4}
\email{\tt colliand@math.toronto.edu}

\author{M. Grillakis}
\address{Department of Mathematics, University of Maryland,
College Park, MD 20742} \email{\tt mng@math.umd.edu}

\author{N. Tzirakis}
\address{Department of Mathematics, University of Toronto, Toronto, Ontario, Canada M5S 2E4}
\email{tzirakis@math.toronto.edu}
\date{13 March 2006}


\begin{abstract}
We prove global well-posedness for low regularity data for the $L^2-critical$ defocusing nonlinear
 Schr\"odinger equation (NLS) in 2d. More precisely we show that a global solution exists for initial data in the Sobolev
 space $H^{s}(\mathbb R^2)$ and any $s>\frac{2}{5}$. This improves the previous result of Fang and Grillakis where global
 well-posedness was established for any $s \geq \frac{1}{2}$. We use the $I$-method to take advantage of the conservation laws
 of the equation. The new ingredient is an interaction Morawetz estimate similar to one that has been used to obtain 
global well-posedness and scattering for the cubic NLS in 3d. The derivation of the estimate in our case is technical since
 the smoothed out version of the solution $Iu$ introduces error terms in the interaction Morawetz inequality. A byproduct of the method
 is that the $H^{s}$ norm of the solution obeys polynomial-in-time bounds.
\end{abstract}
\maketitle

\section{Introduction}

In this paper we study the $L^{2}-$critical Cauchy problem

\begin{equation}\label{ivp1}
\left\{
\begin{matrix}
iu_{t}+ \Delta u -|u|^{2}u=0, & x \in {\mathbb R^2}, & t\in {\mathbb R},\\
u(x,0)=u_{0}(x)\in H^{s}({\mathbb R^2}).
\end{matrix}
\right.
\end{equation}
The problem is known to be locally well-posed for any $s>0$. 
The local well-posedness definition that we use here reads as
follows: for any choice of initial data $u_0 \in H^s$, there exists
a positive time $T = T(\|u_0\|_{H^{s}})$ depending only on the
norm of the initial data, such that a solution to the initial
value problem exists on the time interval $[0,T]$, it is unique in
a certain Banach space of functions $X\subset C([0,T],H^s_{x})$,
and the solution map from $H^s_{x}$ to $C([0,T],H^s_{x})$ depends
continuously on the initial data on the time interval $[0,T]$. If
the time $T$ can be proved to be arbitrarily large, we say that the Cauchy problem is globally well-posed.
A local solution also exists for $L^{2}$ initial data
 but the time of existence depends not only on the $L^2$ norm of the initial data but also on the profile of $u_{0}$. 
For all the above results the reader can look at \cite{tc}, \cite{jb2}, and \cite{tt}. Local in time solutions enjoy mass
conservation
\begin{equation}\label{mass}
\|u(\cdot,t)\|_{L^2(\mathbb R^2)} =
\|u_0(\cdot)\|_{L^2(\mathbb R^2)}.
\end{equation}
Moreover, $H^{1}$ solutions enjoy conservation of the energy
\begin{equation}
E(u)(t)=\frac{1}{2}\int_{\mathbb R^2} |\nabla u(t)|^{2}dx+\frac{1}{4}\int_{\mathbb R^2}
|u(t)|^{4}dx=E(u)(0),
\end{equation} which together with \eqref{mass} and the local theory immediately
yields global-in-time well-posedness for \eqref{ivp1}
with initial data in $H^1.$  
The $L^{2}$ local well-posedness of the IVP \eqref{ivp1} in an interval $[0,T]$ and the conservation
 of the $L^{2}$ norm, cannot be immediately used to 
 give global-in-time solutions as in the case of the finite energy data. As we said in this case $T=T(u_{0})$
 and the lifetime of the local-in-time result can approach zero for
 fixed $L^2$ norm. For the
 focusing case it is known that large mass solutions can blow-up in finite time. Nevertheless in the defocusing case no blowup solutions
 are known and thus it is conjectured that \eqref{ivp1} is globally well-posed and scatters for $L^{2}$ initial data. 
In other words we expect the solution of the nonlinear equation to scatter to a free solution $e^{it\Delta}u_{\pm}$ as 
$t \rightarrow \pm \infty$ for some $u_{\pm} \in L_{x}^{2}(\Bbb R^2)$ in the sense that
$$\lim_{t\rightarrow \pm \infty}\|u(t)-e^{it\Delta}u_{\pm}\|_{L_{x}^{2}(\Bbb R^2)}=0.$$
Conversely, given any $u_{\pm}$ in $ L_{x}^{2}(\Bbb R^2)$ there exists a solution which scatters to it in the sense above, thus giving
 rise to well-defined wave and scattering operators. This conjecture is known to be true in the case that the initial data
 has sufficiently small mass. For details see, \cite{tc}.

For solutions below the energy threshold the first result was
established by J. Bourgain, \cite{jb1}. Bourgain decomposed the initial data into low frequencies and high frequencies and 
estimated seperately the evolution of low and high frequencies. He showed that the solution is globally well-posed with
 initial data in $H^s(\Bbb R^2)$ for any $s>\frac{3}{5}$. Moreover if we denote with $S_{t}$ the nonlinear flow and with 
$S(t)=e^{it\Delta}u_{0}$ the linear group, Bourgain's method shows in addition that $\left( S_t-S(t)\right)u_{0} \in H^{1}(\Bbb R^2)$ 
for all times provided $u_0 \in H^s, s > \frac{3}{5}$. Further
improvements (but without the $H^1$ proximity of the linear and
nonlinear flows) where given in  \cite{ckstt2}
and \cite{fg}, where the authors used the ``$I$-method'' that we describe below. 
Recently T. Tao, M. Visan, and X. Zhang proved \cite{tvz} global well-posedness and scattering
 for the $L^2$-critical problem in all dimensions $n \geq 3$, assuming radially symmetric initial data. They used the reductions in \cite{tvz1}
 to eliminate blow-up solutions that are almost periodic modulo scaling. As in \cite{ckstt5} they obtained a frequency-localized
 Morawetz estimate and exclude a ``mass evacuation scenario'' in order to conclude the argument. Their argument cannot be extended
 easily in low dimensions or without the radial assumption on the initial data. We on the other hand consider the general problem
 in 2d with general initial data but we only relax the regularity
 requirements of the initial data, being unable so far to
 prove the result for initial data in $L^2$. 
\\
\\
We use the $I$-method and we follow closely the argument in \cite{fg}
(see also \cite{ckstt2}, \cite{ckstt4}) which is based on the almost conservation
of a certain modified energy functional. 
The idea is to replace
the conserved quantity $E(u)$ which is no longer available for
$s<1$, with an ``almost conserved'' variant $E(Iu)$ where $I$ is a
smoothing operator of order $1-s$ which behaves like the identity
for low frequencies and like a fractional integral operator for
high frequencies. Thus the operator $I$ maps $H_{x}^{s}$ to
$H_{x}^{1}$. Notice that $Iu$ is not a solution to \eqref{ivp1}
and hence we expect an energy increment. This increment is in fact
quantifying $E(Iu)$ as an ``almost conserved'' energy. The key is
to prove that on intervals of fixed length, where local
well-posedness is satisfied, the increment of the modified energy
$E(Iu)$ decays with respect to a large parameter $N$. (For the
precise definition of $I$ and $N$ we refer the reader to Section
$2$.) This requires delicate estimates on the commutator between
$I$ and the nonlinearity. In dimensions 1 and 2, where the
nonlinearity is algebraic, one can write the commutator explicitly
using the Fourier transform, and control it by multi-linear
analysis and bilinear estimates. The statement of the our main result is:

\begin{thm}\label{main2d}
The initial value problem \eqref{ivp1} is globally
well-posed in $H^{s}(\mathbb R^2)$, for any $1>s>\frac{2}{5}$. Moreover the
solution satisfies
$$\sup_{[0,T]}\|u(t)\|_{H^{s}(\mathbb R^2)} \leq C(1+T)^{\frac{3s(1-s)}{2(5s-2)}}$$
where the constant $C$ depends only on the index $s$ and $\|u_{0}\|_{L^{2}}$.
\end{thm}

\begin{rem}
  We view this result as another incremental step towards the
  conjecture that $L^2$ initial data $u_0$ evolves under \eqref{ivp1}
  to a global-in-time solution with $\| u \|_{L^4{(\R_t \times \R^2_x)}}
  < C(u_0)$. Recent work \cite{CKSTT:Angles}, based upon the second
  modified energy and certain angular refinements of the bilinear
  Strichartz estimate, has improved the energy increment
  quantification (see \eqref{EnergyIncrement}) from $N^{-3/2+}$ to
$N^{-2+}$. This improvement is applied in \cite{CKSTT:Angles}
following the globalizing scheme in \cite{ckstt2} (which does not rely
upon any monotonicity or Morawetz-type inputs) to prove that
\eqref{ivp1} is globally well-posed in $H^s (\R^2)$ for $s >
\frac{1}{2}$. Similarities between the proofs of the almost
conservation estimate \eqref{EnergyIncrement} and the almost Morawetz
increment estimate \eqref{AlmostMorawetz} suggest that angle
refinements as in \cite{CKSTT:Angles} could possibly improve the
decay in \eqref{AlmostMorawetz} below $N^{-1+}$, possibly to
$N^{-3/2+}$. If true, this would improve the global well-posedness
range to $s > \frac{4}{13}$.
\end{rem}

The basic ingredient in our proof is an \textit{a priori}
interaction Morawetz-type estimate for the ``approximate solution'' $Iu$ to the initial value problem
\begin{equation}\label{aivp1}
\left\{
\begin{matrix}
iIu_{t}+ \Delta Iu -I(|u|^{2}u)=0 & x \in {\mathbb R^2},& t\in {\mathbb R}\\
Iu(x,0)=Iu_{0}(x)\in H^{1}({\mathbb R^2}).
\end{matrix}
\right.
\end{equation}
For the original system \eqref{ivp1} it has been shown in \cite{fg} that solutions satisfy the following a priori bound
$$\|u\|_{L_{T}^{4}L_{x}^{4}}^{4} \lesssim T^{\frac{1}{2}}\sup_{t\in [0,T]}\|u(t)\|_{\dot{H}^{\frac{1}{2}}}^2\|u_0\|_{L^{2}}^{2}.$$
This estimate follows from a ``two-particle'' Morawetz inequality. ``Two-particle'' Morawetz estimates first appeared in three
 dimensions in \cite{ckstt4}. 
Roughly speaking, we refer to Morawetz inequalities as 
monotonicity formulae that take advantage of the conservation of momentum
$$\vec{p}(t)=\int_{\Bbb R^{n}}\Im(\bar{u}(x)\nabla u(x))dx.$$
In dimensions $n\geq 3$ the classical Morawetz inequality relies ultimately on the fact
that the tempered distribution $\Delta \Delta |x|$ is well-defined and nonnegative. In particular for $n=3$ we have that
$\Delta \Delta |x|= \delta_{0}$. The case $n=2$ is more delicate and this is the novelty of the approach in \cite{fg}.
In this paper we improve the above inequality to
$$\|u\|_{L_{T}^{4}L_{x}^{4}}^{4} \lesssim T^{\frac{1}{3}}\sup_{t\in
  [0,T]}\|u(t)\|_{\dot{H}^{\frac{1}{2}}}^2\|u_0\|_{L^{2}}^{2}
+CT^{\frac{1}{3}}\|u_0\|_{L^{2}}^{4}$$
where $C$ is a large constant.
For initial data below $H^{\frac{1}{2}}$ both of these estimates are not useful. This is mainly the limitation of the result
 in \cite{fg}. To avoid this difficulty we can introduce
 the $I$-operator and hope to get an apriori estimate of the form
\begin{equation}
\label{IuMorawetzThird}\|Iu\|_{L_{T}^{4}L_{x}^{4}}^{4} \lesssim
T^{\frac{1}{3}}\sup_{t\in
  [0,T]}\|Iu(t)\|_{\dot{H}^{1}}^{2}\|Iu_0\|_{L^{2}}^{2}
+CT^{\frac{1}{3}}\|u_0\|_{L^{2}}^{4} + Error.
\end{equation}
Then the restriction $s \geq \frac{1}{2}$ is not present and, in
principle, one can improve the result in \cite{fg}.
Of course, we have  to show that the $Error$ terms are negligible in some
sense. More precisely we show that on the local well-posedness time
interval the error terms are very small. The proof of this fact relies
on multilinear harmonic analysis
 estimates (similar to those used to prove almost conservation
 bounds in \cite{ckstt4}) and is given in Section 4. Before we outline the general method 
of the paper we define the following Banach space:
$$\label{ZI} Z_I(J) :=\sup_{(q,r) \ \ admissible} \|\langle \nabla \rangle
Iu\|_{L^q_tL^r_x(J\times \mathbb R^2)}$$
where a pair $(q,r)$ is said to be admissible if
$\frac{1}{q}+\frac{1}{r}=\frac{1}{2}$, $2 \leq r,q \leq \infty$. We
will sometimes write $Z_I (t) $ to denote $Z_I([0,t])$.
\\
\\
Now fix a large value of time $T_{0}$. If $u$
is a solution to \eqref{ivp1} on the time interval $[0,T_{0}],$
then
$u^\lambda(x)=\frac{1}{\lambda}u(\frac{x}{\lambda},\frac{t}{\lambda^{2}})$
is a solution to the same equation on $[0, \lambda^{2}T_{0}]$. We
choose the parameter $\lambda>0$ so that
$E(Iu_{0}^{\lambda})=O(1)$. Using Strichartz estimates we show that if $J=[a,b]$ and
$\|Iu^{\lambda}\|_{L^{4}_tL^{4}_x(J\times
\Bbb R^2 )}^{4}<\mu$, where $\mu$ is a small universal
constant, then
$$\sup_{(q,r) \ \ admissible} \|\langle \nabla \rangle
Iu^{\lambda}\|_{L^q_tL^r_x(J\times \mathbb R^2)} \lesssim \|I
u^\lambda(a)\|_{H^1} \lesssim 1.$$ 
Moreover in  this same time interval where
the problem is well-posed, we can prove the ``almost
conservation law'' 
\begin{equation}\label{acl}
    |E(Iu^{\lambda})(b) - E(Iu^{\lambda})(a)| \lesssim N^{-\frac{3}{2}}\|Iu^\lambda(a)\|_{H^1}^{4}+N^{-2}\|Iu^\lambda(a)\|_{H^1}^{6}
\lesssim N^{-\frac{3}{2}}.
\end{equation} 
For the arbitrarily large
interval $[0,\lambda^{2}T_{0}] $ we do not have that
$$\|Iu^{\lambda}\|_{L^{4}_tL^{4}_x([0,\lambda^{2}T_{0}]\times
\mathbb R^2 )}^{4}<\mu.$$
But we can partition the  arbitrarily large interval
$[0,\lambda^{2}T_{0}]$ into $L$ intervals where the local theory
uniformly applies. $L=L(N,T)$  is finite and defines  the number of the intervals in the partition
that will make the Strichartz
$L^{4}_tL^{4}_x$ norm of $Iu$ less than
$\mu$.  Since $E(Iu^{\lambda})$ controls the $H^1$ norm of $Iu$ we have by (\ref{acl}) that
$$\|Iu^{\lambda}\|_{H^{1}}^2 \lesssim LN^{-\frac{3}{2}}.$$
To maintain the bound $\|Iu^{\lambda}\|_{H^1} \lesssim 1$ we choose
$$L(N,T)\sim N^{\frac{3}{2}}$$ and this condition will ultimately lead
to the requirement $s>\frac{2}{5}$. For the detailed proof, see Section 5.

The rest of the paper is organized as follows. In Section 2 we
introduce some notation and state important propositions that we
will use throughout the paper. In Section 3 we
prove the local well-posedness theory for $Iu,$ and the main
estimates that we use to prove the decay of the increment of the
modified energy. The decay itself is obtained in Section 3. In Section 4 we prove the ``almost Morawetz'' inequality which
 is the heart of our argument. Finally in Section 5 we give the details of the proof of global
well-posedness stated in Theorem \ref{main2d}.
\\
{ \bf Acknowledgments.} We would like to thank D. De Silva, Y. Fang,
M. Keel, N. Pavlovic, G. Staffilani, H. Takaoka and T. Tao 
for discussions related to this work.

\section{Notation}

In what follows we use $A \lesssim B$ to denote an estimate of the form $A\leq CB$ for some constant $C$.
If $A \lesssim B$ and $B \lesssim A$ we say that $A \sim B$. We write $A \ll B$ to denote
an estimate of the form $A \leq cB$ for some small constant $c>0$. In addition $\langle a \rangle:=1+|a|$ and
$a\pm:=a\pm \epsilon$ with $0 < \epsilon \ll 1$.
\\
\\
We use $\mathcal{S}$ to denote the Schwartz class. 
We use $L^r_x(\mathbb R^2)$ to denote the Lebesgue space of functions $f :
\R^2 \rightarrow \C$ whose norm $$\|f\|_{L^r_x}:=\left(
\int_{\R^2}|f(x)|^r dx \right)^{\frac{1}{r}}$$ is finite, with the
usual modification in the case $r=\infty.$ We also define the
space-time spaces $L^q_tL^r_x$ by $$\|u\|_{L^q_{t\in J}L^r_x} :=
\left(\int_J \|u\|_{L^r_x}^q dt\right)^{\frac{1}{q}}$$ for any
space-time slab $J \times \mathbb R^2,$ with the usual modification when
either $q$ or $r$ are infinity. When $q=r$ we abbreviate
$L^q_tL^r_x$ by $L^q_{t,x}.$ Sometimes we will write $L^q_T$ to denote
$L^q_{t \in [0,T]}$.

\begin{defn} A pair of exponents $(q,r)$ is called admissible  in
    $\mathbb R^2$ if
$$ \frac{1}{q}+\frac{1}{r} = \frac{1}{2}, \ \ \ 2 \leq q,r \leq \infty.$$
\end{defn}
We recall the Strichartz estimates \cite{gv}, \cite{kt} (and the
references therein).
\begin{prop}Let $(q,r)$ and $(\tilde{q},\tilde{r})$ be
any two admissible pairs. Suppose that $u$ is a solution to
\begin{equation*}
\left\{
\begin{matrix}
i\partial_t u+ \Delta u -G(x,t)=0, & (t,x) \in J \times \mathbb R^2 \\
u(x,0)=u_{0}(x).\\
\end{matrix}
\right.
\end{equation*} 
Then we have the estimate \begin{equation}\label{S}
\|u\|_{L^q_tL^r_x (J\times \mathbb R^2)} \lesssim
\|u_0\|_{L^2(\mathbb R^2)}+\|G\|_{L^{\tilde q^{\prime}}_tL^{\tilde r^{\prime}}_x
(J\times\mathbb R^2)}
\end{equation}
with the prime exponents denoting H\"older dual exponents.
\end{prop}
We define the Fourier transform of $f(x) \in L_{x}^{1}$ by
$$\hat{f}(\xi)=\int_{\mathbb R^2} e^{-2\pi i\xi x}f(x)dx.$$
For an appropriate class of functions the following Fourier inversion formula holds:
$$f(x)=\int_{\mathbb R^2} e^{2\pi i\xi x}\hat{f}(\xi)(d\xi).$$
Moreover we know that the following identities are true:\\
\begin{enumerate}
\item $\|f\|_{L^{2}}=\|\hat{f}\|_{L^{2}}$, (Plancherel)\\
\item $\int_{\mathbb R^2} f(x)\bar{g}(x)dx=\int_{\mathbb R^2} \hat{f}(\xi)\bar{\hat{g}}(\xi) (d\xi)$, (Parseval)\\
\item $\widehat{fg}(\xi)=\hat{f}\star \hat{g}(\xi)=\int_{\mathbb R^2} \hat f(\xi-\xi_{1})\hat g(\xi_{1})d\xi_{1}$, (Convolution).\\
\end{enumerate}
We also
define the fractional differentiation operator $D^\alpha$
for any real $\alpha$ by $$\widehat{D^\alpha u}(\xi):=
|\xi|^\alpha \hat{u}(\xi)$$ and analogously
$$\widehat{\langle D \rangle^\alpha u}(\xi):=
\langle\xi \rangle^\alpha \hat{u}(\xi).$$ 
We then define the
inhomogeneous Sobolev space $H^s$ and the homogeneous Sobolev
space $\dot{H}^s$ by $$\|u\|_{H^s} = \|\langle D \rangle^s
u\|_{L^2_x}; \ \ \ \ \|u\|_{\dot{H}^s} =\|D^s u\|_{L^2_x}.$$

\section{The I-method and the proof of Theorem \ref{main2d}}
We shall also need some Littlewood-Paley theory, \cite{es}. The reader must have in mind that wherever
 in this paper we restrict the functions in frequency we do it in a smooth way using the Littlewood-Paley projections. In particular, let
$\eta(\xi)$ be a smooth bump function supported in the ball $|\xi|
\leq 2,$ which is equal to one on the unit ball. Then, for each
dyadic number $M$ we define the Littlewood-Paley operators

\begin{align*}\widehat{P_{\leq M} f}(\xi) &= \eta(\xi/M) \hat{f}(\xi),\\
\widehat{P_{> M} f}(\xi) &= (1-\eta(\xi/M)) \hat{f}(\xi),\\
\widehat{P_{M} f}(\xi) &= (\eta(\xi/M)-\eta(2\xi/M)) \hat{f}(\xi).\\
\end{align*}
Similarly, we can define $P_{<M}, P_{\geq M}.$  The Littlewood-Paley
decomposition we write, at least formally, is $u=\sum_{M}P_{M}u$. 
For convenience we abbreviate the Littlewood-Paley operator $P_{M}$ by
 $u_{M}$ or even $u_{j}$ when its meaning is clear from the context. We can write
 $u=\sum u_{j}$ and obtain bounds on each piece seperately or by examining the
 interactions of the several pieces. We can recover information for the original function $u$ by applying the
 Cauchy-Schwartz inequality and using the Littlewood-Paley Theorem or the cheap Littlewood-Paley inequality
$$\|P_{N}u\|_{L^p} \lesssim \|u\|_{L^{p}}$$
for any $1\leq p \leq \infty$. 
Since this process is fairly standard
 we will omit the details of the argument throughout the paper.
We also recall the
definition of the operator $I$. For $s<1$ and a parameter $N \gg 1$
let $m(\xi)$ be the following smooth
monotone multiplier:
\[m(\xi):= \left\{\begin{array}{ll}
1 & \mbox{if $|\xi|<N$}\\
(\frac{|\xi|}{N})^{s-1} & \mbox{if $|\xi|>2N$}.
\end{array}
\right.\] We define the multiplier operator $I:H^{s} \rightarrow
H^{1}$ by
$$\widehat{Iu}(\xi)=m(k)\hat{u}(\xi).$$
The operator $I$ is smoothing of order $1-s$ and we have that:
\begin{equation}
\|u\|_{H^{s_{0}}} \lesssim \|Iu\|_{H^{s_{0}+1-s}} \lesssim N^{1-s}\|u\|_{H^{s_{0}}}
\end{equation}
for any $s_{0} \in {\mathbb R}$.
\\
\\
We set
\begin{equation}\label{first}E^{1}(u)=E(Iu),\end{equation} where $$E(u)(t)=\frac{1}{2}\int
|u_{x}(t)|^{2}dx+\frac{1}{6}\int |u(t)|^{6}dx=E(u_{0}).$$ We refer
to $E^{1}(u)$ as the first modified energy. 
\subsection{Modified Local Well-Posedness}
\begin{prop}\label{lwp}
Define the quantity $\mu([0,T])=\int_{0}^{T}\int_{\mathbb R^2}|Iu(x,t)|^4dxdt.$ If $\mu([0,T]) < \mu_{0}$ where $\mu_{0}$
 is a universal constant then for any $s>0$ the initial value problem \eqref{aivp1} is locally
 well-posed and we have that
$$ Z_{I}([0,T]) :=\sup_{(q,r) \ \ admissible} \|\langle D \rangle
Iu\|_{L^q_tL^r_x([0,T]\times \mathbb R^2)} \lesssim \|\langle D \rangle Iu_{0}\|_{L^{2}}.$$
\end{prop}

\begin{proof}
By standard well-posedness theory, see for example \cite{fg}, it is enough to show 
$$ Z_{I}([0,T]) :=\sup_{(q,r) \ \ admissible} \|\langle D \rangle
Iu\|_{L^q_tL^r_x([0,T]\times \mathbb R^2)} \lesssim \|\langle D \rangle Iu_{0}\|_{L^{2}}.$$
By the Duhamel formula we have that an equivalent representation for the solution of \eqref{aivp1} and is given by
$$Iu(x,t)=e^{it\Delta}Iu_{0}-i\int_{0}^{T}e^{i(t-s)\Delta}I(|u|^4u)(s)ds.$$ Applying the $\langle D \rangle$ operator in
the equation \eqref{aivp1} and using the Strichartz estimates \eqref{S} we have that
$$Z_{I}\lesssim \|Iu_{0}\|_{H^{1}}+\|\langle D \rangle I(|u|^2u)\|_{L_{t}^{\frac{4}{3}}L_{t}^{\frac{4}{3}}}
\lesssim \|Iu_{0}\|_{H^{1}}+\|\langle D \rangle Iu\|_{L_{t}^{4}L_{t}^{4}}\|u\|_{L_{t}^{4}L_{t}^{4}}^{2}.$$
Note that we have used Leibnitz's rule for fractional derivatives in the previous step. Indeed the multiplier $\langle D \rangle I$
 is increasing for any $s \geq 0$. Using this fact one can modify the proof of the usual Leibnitz rule for fractional derivatives
 and prove it also for $\langle D \rangle I$. Thus,
\begin{equation}\label{Z} 
Z_{I}\lesssim \|Iu_{0}\|_{H^{1}}+Z_{I}\|u\|_{L_{t}^{4}L_{t}^{4}}^{2}.
\end{equation}
Now recall the definition of $I$. We write $u=u_{<N}+\sum_{j=1}^{\infty}u_{k_j}$, where $u_{<N}$ has spatial frequency support
 on $\langle \xi \rangle \leq N$ and the $u_{k_j}$ have support $\langle
 \xi_{j} \rangle \sim N_{j}=2^{k_{j}}$ where $k_j$ are consecutive
 integers starting with $[\log N]$ indexed by
  $j=1,2,...,$. Note that
$$
\|u_{k_j}\|_{L_{t}^{4}L_{t}^{4}} \lesssim
N_{j}^{1-s}N^{s-1}\|Iu_{k_j}\|_{L_{t}^{4}L_{t}^{4}}, ~ \mbox{if $j=1,2,...$}
$$
By the triangle inequality we have that
$$\|u\|_{L_{t}^{4}L_{t}^{4}} \leq \|u_{<N}\|_{L_{t}^{4}L_{t}^{4}}+\sum_{j=1}^{\infty}\|u_{k_j}\|_{L_{t}^{4}L_{t}^{4}}
=\|Iu_{<N}\|_{L_{t}^{4}L_{t}^{4}}+\sum_{j=1}^{\infty}\|u_{k_j}\|_{L_{t}^{4}L_{t}^{4}}^{\epsilon}
\|u_{k_j}\|_{L_{t}^{4}L_{t}^{4}}^{1-\epsilon}.$$
By the definition of the $u_{j}$'s the following estimates are true
\begin{equation}
\|Iu_{<N}\|_{L_{t}^{4}L_{t}^{4}} \lesssim \|Iu\|_{L_{t}^{4}L_{t}^{4}}
\end{equation} 
\begin{equation}
\|u_{k_j}\|_{L_{t}^{4}L_{t}^{4}} \lesssim N_{j}^{1-s}N^{s-1}\|Iu_{k_j}\|_{L_{t}^{4}L_{t}^{4}},\ \ j=1,2,...
\end{equation}
\begin{equation}
\|u_{k_j}\|_{L_{t}^{4}L_{t}^{4}} \lesssim N_{j}^{-s}N^{s-1}\|\langle D \rangle Iu_{k_j}\|_{L_{t}^{4}L_{t}^{4}},\ \ j=1,2,....
\end{equation} 
Combining all these estimates we get that
$$\|u\|_{L_{t}^{4}L_{t}^{4}} \lesssim \|Iu\|_{L_{t}^{4}L_{t}^{4}}+\sum_{j=1}^{\infty}N_{j}^{-s+\epsilon}N^{s-1}
\|\langle D \rangle Iu_{k_j}\|_{L_{t}^{4}L_{t}^{4}}^{1-\epsilon}\|Iu_{k_j}\|_{L_{t}^{4}L_{t}^{4}}^{\epsilon}$$
and after eliminating $N^{s-1} \leq 1$ we have
$$\|u\|_{L_{t}^{4}L_{t}^{4}} \lesssim \|Iu\|_{L_{t}^{4}L_{t}^{4}}+\sum_{j=1}^{\infty}N_{j}^{-s+\epsilon}
\|\langle D \rangle Iu_{k_j}\|_{L_{t}^{4}L_{t}^{4}}^{1-\epsilon}\|Iu_{k_j}\|_{L_{t}^{4}L_{t}^{4}}^{\epsilon}.$$
Now if we apply the cheap Littlewood-Paley inequality 
$$\|u_{k_j}\|_{L^p} \lesssim \|u\|_{L^p}$$
for any $1 \leq p \leq \infty$ and sum,
we have that for any $s>\epsilon$ we get that
$$\|u\|_{L_{t}^{4}L_{t}^{4}} \lesssim \|Iu\|_{L_{t}^{4}L_{t}^{4}}+\|\langle D \rangle Iu\|_{L_{t}^{4}L_{t}^{4}}^{1-\epsilon}
\|Iu\|_{L_{t}^{4}L_{t}^{4}}^{\epsilon}.$$
Putting all these into equation \eqref{Z} we have
$$Z_{I} \lesssim \|Iu_{0}\|_{H^{1}}+Z_{I}\|Iu\|_{L_{t}^{4}L_{t}^{4}}^2+Z_{I}^{3-2\epsilon}\|Iu\|_{L_{t}^{4}L_{t}^{4}}^{2\epsilon}$$
Now if we pick a time $T$ such that $\mu([0,T]) < \mu_{0} \ll 1$ we have that
$$ Z_{I}([0,T]) :=\sup_{(q,r) \ \ admissible} \|\langle D \rangle
Iu\|_{L^q_tL^r_x([0,T]\times \mathbb R^2)} \lesssim \|\langle D \rangle Iu_{0}\|_{L^{2}}.$$
\end{proof}
We finish this section with the almost conservation of the modified energy.
\begin{prop} \label{ac}
If $H^s \ni u_0 \longmapsto u(t)$ with $\frac{1}{2}>s>\frac{1}{3}$ solves \eqref{ivp1} for all $t \in [0,T_{lwp}]$ 
where $T_{lwp}$ is the time that Proposition \ref{lwp} applies. Then
\begin{equation}
 \sup_{t \in [0,T_{lwp}]} |E[ I_N u(t) ]| \leq |E[I_N u(0) ]|
  + C N^{-\frac{3}{2}-} \| I_N \langle D \rangle u(0) \|_{L^2_x}^4 + C N^{-2-}
  \|I_N \langle D \rangle u(0) \|_{L^2_x}^6.
\end{equation}
In particular when $\|I_N \langle D \rangle u(0) \|_{L^2_x} \lesssim 1$ we have that
\begin{equation}
\label{EnergyIncrement}
\sup_{t \in [0,T_{lwp}]} |E[ I_N u(t) ]| \leq
 |E[I_N u(0) ]|+C N^{-\frac{3}{2}-} \lesssim  N^{-\frac{3}{2}-}.
\end{equation}
\end{prop}
\begin{proof}
This is Proposition 3.7 in \cite{crsw}. The restriction for
$s>\frac{1}{3}$ appears in Case 3 in the proof.
\end{proof}

\section{The almost Morawetz estimate.}
For what follows we sometimes abbreviate $u_{i}=u(x_{i})$ where $u_{i}$ is a solution to 
\begin{equation}
iu_{t}+\Delta u=|u|^2u, ~ (x_{i},t) \in \R^2 \times [0,T].
\end{equation} 
Here $x_i
\in \R^n$, not a coordinate. 
In this section we wish to prove the {'almost Morawetz'} estimate. For this consider 
$a: \mathbb R^n \rightarrow \Bbb R$, a convex and locally
integrable function of polynomial growth.

\begin{thm}\label{main}Let $u \in L^\infty_{[0,T]} \mathcal{S}_x$ be a solution to the NLS
\begin{equation}\label{NLSq}
iu_{t}+\Delta u=|u|^2u, ~ (x,t) \in \R^2 \times [0,T]
\end{equation} 
and 
$Iu \in L^\infty_{[0,T]} \mathcal{S}_x$ be a solution to the $I$-NLS
\begin{equation}
iIu_{t}+\Delta Iu=I(|u|^2u), ~ (x,t) \in \R^2 \times [0,T].
\end{equation} 
Then,
\begin{align}\label{amorerror}
\|Iu\|_{L_{T}^{4}L_{x}^{4}}^{4} &\lesssim
T^{\frac{1}{3}}\sup_{[0,T]}\|Iu\|_{\dot{H}^{1}}^1 \|Iu\|_{L^{2}}^{3}+T^{\frac{1}{3}}\|u_{0}\|_{L^{2}}^{4}+\\ \nonumber & 
T^{\frac{1}{3}}\int_0^T\int_{\mathbb R^2 \times \mathbb R^2}\nabla a \cdot
\{\widetilde{\CN}_{bad},Iu(x_{1},t)Iu(x_{2},t)\}_{p}dx_1dx_2dt.
\end{align}
where 
$$\widetilde{\CN}_{bad}=\sum_{i=1}^{2}\left(I(|u_{i}|^2u_{i})-|Iu_{i}|^2Iu_{i}\right)\prod_{j=1,j \ne
i}^{2}Iu_{j}$$ and $\{\cdot\}_{p}$ is the momentum bracket defined by 
$$\{f,g\}_p=\Re(f\overline{\nabla g} -
g\overline{\nabla f}).$$
In particular, on a time interval $J_{k}$ where the local well-posedness Proposition \ref{lwp} holds we have that
\begin{equation}
\label{AlmostMorawetz}
 \int_{J_k}\int_{\mathbb R^2 \times \mathbb R^2}\nabla a \cdot
\{N_{bad},Iu(x_{1},t)Iu(x_{2},t)\}_{p}dx_1dx_2dt \lesssim \frac{1}{N^{1-}}Z_I^{6}(J_k).
\end{equation}
\end{thm}
Toward this aim, we recall the idea of the proof of the
interaction Morawetz estimate for the defocusing nonlinear
 cubic Schr\"odinger equation in three space dimensions
 \cite{ckstt4}. 
We present the
 result using a different argument involving a tensor of Schr\"odinger
 solutions that emerged from a conversation
 between Andrew Hassell and Terry Tao. 
We will establish the ``almost Morawetz'' estimate that we need in this paper, along the lines
 of this new point of view. In all of our arguments we will assume smooth solutions. This will simplify
 the calculations and will enable us justify the steps in the subsequent proofs. The local well-posedness theory
 and the perturbation theory \cite{tc} that has been established for this problem can be then applied to approximate
 the $H^{s}$ solutions by smooth solutions and conclude the proofs. For most of the calculations in this section the reader can 
consult \cite{ckstt5}, \cite{tt}.  

Let start with a solution to the NLS
\begin{equation}\label{NLS}
iu_{t}+\Delta u=\widetilde{\CN}(u), \ \ \ (x,t) \in \mathbb R^n \times [0,T]
\end{equation}
with $u$ Schwartz-class-in-space and $\widetilde{\CN}$
such that there exist a defocusing potential $G$, (meaning $G$
positive) such that $$\{\widetilde{\CN},u\}_{p}^{j}=-\partial_{j}G.$$ 
Let's define also the momentum density
$$T_{0j}=2\Im (\bar{u}\partial_{j}u)$$
for $j=1,2,...,n$, and the linearized momentum current
$$L_{jk}=-\partial_{j}\partial_{k}(|u|^2)+4\Re (\overline{\partial_{j}u}\partial_{k}u).$$ 
A computation shows that
$$\partial_{t}T_{0j}+\partial_{k}L_{jk}=2\{\widetilde{\CN},u\}_{p}^{j}$$
where we have adopted Einstein's summation convention. Notice also that in our case where $\widetilde{\CN}=|u|^2u$ we have that 
$\{\widetilde{\CN},u\}_{p}^{j}=-\partial_{j}G$, 
where $G=\frac{1}{2}|u|^{4}$. By integrating in space we have that the total momentum
 is conserved in time,
$$\int_{\mathbb R^{n}}T_{0j}(x,t)dx=C.$$
We recall the {\it{generalized virial identity}} \cite{LinStrauss}.
\begin{prop}\label{wmor}
If $a$ is convex and $u$ is 
a smooth solution to equation \eqref{NLS} on $[0,T]\times \mathbb R^{n}$ with a defocusing potential $G$. 
Then, the following inequality holds:
\begin{equation}\label{Mor}
\int_0^T\int_{\Bbb R^n}(-\Delta \Delta a)|u(x,t)|^2dxdt \lesssim 
\sup_{[0,T]}|M_{a}(t)|,
\end{equation}
where $M_{a}(t)$ is the Morawetz action and is given by
\begin{equation}
\label{Mat}M_a(t)= 2 \int_{\Bbb R^n}\nabla a(x)
\cdot \Im(\overline{u}(x)\nabla u)dx.
\end{equation}
\end{prop}
\begin{proof}
We can write the Morawetz action as
$$M_{a}(t)=\int_{\Bbb R^n} \partial_{j}aT_{0j}.$$
Then
$$\partial_{t}M_a(t)=\int_{\Bbb R^n} \partial_{j}a\partial_{t}T_{0j}=\int_{\Bbb R^n} 
\partial_{j}a\left( -\partial_{k}L_{jk}+2\{\widetilde{\CN},u\}_{p}^{j}\right)=$$
$$\int_{\Bbb R^n} 
\partial_{j}a\left( -\partial_{k}L_{jk}-2\partial_{j}G\right)=\int_{\Bbb R^n}(\partial_{j}\partial_{k}a)L_{jk}dx
+2\int_{\Bbb R^n}\Delta a Gdx$$
where in the last equality we used integration by parts.
By the definition of $L_{jk}$ we have that
$$\partial_{t}M_a(t)=\int_{\mathbb R^n}(\partial_{j}\partial_{k}a)\partial_{j}\partial_{k}(|u|^2)dx+
4\int_{\mathbb R^n}(\partial_{j}\partial_{k}a)\Re
\left( \partial_{j}\bar{u}\partial_{k}u\right)dx+2\int_{\mathbb
  R^n}\Delta a Gdx.$$
Performing the summations, we record the {\it{generalized virial identity}}
$$\partial_t M_a (t) =-\int_{\Bbb R^n}(\Delta \Delta a) |u|^{2}dx+2\int_{\Bbb R^n}\Delta a Gdx+
4\int_{\mathbb R^n}(\partial_{j}\partial_{k}a)\Re
\left( \partial_{j}\bar{u}\partial_{k}u\right)dx.$$
But since $a$ is convex we have that 
$$4(\partial_{j}\partial_{k}a)\Re \left( \partial_{j}\bar{u}\partial_{k}u\right) \geq 0$$
and the trace of the Hessian of $\partial_{j}\partial_{k}a$ which is $\Delta a$ is positive. Thus,
$$-\int_{\Bbb R^n}(\Delta \Delta a) |u|^{2}dx \leq \partial_{t}M_a(t)$$
and by the fundamental theorem of calculus we have that
\begin{equation}
\int_0^T\int_{\Bbb R^n}(-\Delta \Delta a)|u(x,t)|^2dxdt \lesssim 
\sup_{[0,T]}|M_{a}(t)|. \label{mora}
\end{equation}
\end{proof}
In the case of a solution to an equation with a nonlinearity which
is not associated to a defocusing potential, we immediately obtain
the following corollary.

\begin{cor} \label{cor}
Let $a:\R^d \rightarrow \R$ be convex and $u$ be a smooth solution
to the equation
\begin{equation}\label{NLS2} iu_{t}+\Delta u=\widetilde{\CN}, \ \
\ (x,t) \in \mathbb R^d \times [0,T].
\end{equation}Then, the following inequality
holds
\begin{equation}\label{MorN}
\int_{0}^{T}\int_{\Bbb R^n}(-\Delta \Delta a)|u(x,t)|^2dxdt +2
\int_0^T\int_{\R^d} \nabla a \cdot \{\widetilde{\CN},u\}_pdxdt
\lesssim |M_{a}(T) - M_a (t)|,
\end{equation} where $M_{a}(t)$ is the Morawetz action
corresponding to $u$.
\end{cor}

\subsection{Interaction Morawetz inequality in three dimensions.}
Now we consider the interaction Morawetz inequality. Let $u_{i}$, $\widetilde{\CN}_{i}$ 
be solutions to (\ref{NLS}) in $n_{i}-$spatial dimensions
 and suppose we have as before momentum conservation with a defocusing potential. Define the tensor product 
$u:=(u_{1}\otimes u_{2})(t,x)$ for $x$ in 
$$\mathbb R^{n_{1}+n_{2}}=\{(x_{1},x_{2}): x_{1} \in \mathbb R^{n_{1}}, x_{2} \in \mathbb R^{n_{2}}\}$$ 
by the formula 
$$(u_{1}\otimes u_{2})(t,x)=u_{1}(x_{1},t)u_{2}(x_{2},t).$$
It can be easily verified that if $u_1$ solves (\ref{NLS}) with forcing term $\widetilde{\CN}_1$ and  $u_2$ solves (\ref{NLS}) 
with forcing term $\widetilde{\CN}_2$, 
then $u_1 \otimes u_2$ solves (\ref{NLS}) with forcing term $\widetilde{\CN}=\widetilde{\CN}_{1}\otimes u_2+
\widetilde{\CN}_{2}\otimes u_1$. Since
$$\{\widetilde{\CN}_{1}\otimes u_2+\widetilde{\CN}_{2}\otimes u_1,u_1 \otimes u_2\}_{p}=
\left(\{\widetilde{\CN}_1,u_1\}_{p}\otimes |u_2|^2, 
\{\widetilde{\CN}_2,u_2\}_{p}\otimes |u_1|^2\right)$$
we have the important fact that the tensor product of the defocusing semilinear Schr\"odinger equation is also defocusing in
the sense that
$$\{\widetilde{\CN}_{1}\otimes u_2+\widetilde{\CN}_{2}\otimes u_1,u_1 \otimes u_2\}_{p}=-\nabla G$$ 
where $G=G_{1}\otimes |u_2|^2+G_{1}\otimes |u_1|^2$ and $\nabla= (\nabla_{x_1},\nabla_{x_2})$. Thus $G \geq 0$.
Since $u_1 \otimes u_2$ solves (\ref{NLS}) and obeys momentum
conservation with a defocusing potential, we can apply the Proposition
\ref{wmor} and obtain for a convex functions $a$ that
\begin{equation}\label{interact}
\int_0^T\int_{\Bbb R^{n_{1}} \otimes \Bbb R^{n_{2}}}(-\Delta \Delta a)|u_1\otimes u_2|^2(x,t)dxdt \lesssim 
\sup_{[0,T]}|M_{a}^{\otimes_{2}}(t)|
\end{equation}
where $\Delta= \Delta_{x_1}+\Delta_{x_2}$ the Laplacian in $\R^{n_1 +
  n_2}$ and $M_{a}^{\otimes_{2}}(t)$ is the Morawetz
 action that corresponds to $u_1\otimes u_2$ and thus
$$M_{a}^{\otimes_{2}}(t)= 2 \int_{\Bbb R^{n_{1}}\otimes \Bbb R^{n_{2}}}\nabla a(x)
\cdot \Im\left(\overline{u_1\otimes u_2}(x)\nabla (u_1\otimes u_2(x))\right )dx.$$
Now we pick $a(x)=a(x_1,x_2)=|x_1-x_2|$ where $(x_{1},x_2) \in \Bbb R^{3}\times \Bbb R^{3}$. 
Then an easy calculation shows that $-\Delta \Delta a(x_1,x_2)=C\delta(x_1-x_2)$. 
Applying equation \eqref{interact} with this choice of $a$ and choosing $u_{1}=u_{2}$ we get that
$$\int_{0}^{T}\int_{\Bbb R^3}|u(x,t)|^{4}dx \lesssim \sup_{[0,T]}|M_{a}^{\otimes_{2}}(t)|.$$
It can be shown using Hardy's inequality (for details see \cite{ckstt4}) that in 3d
$$\sup_{[0,T]}|M_{a}^{\otimes_{2}}(t)| \lesssim \sup_{[0,T]}\|u(t)\|_{\dot{H}^{\frac{1}{2}}}^2\|u(t)\|_{L^{2}}^2$$
and thus
$$\int_{0}^{T}\int_{\Bbb R^3}|u(x,t)|^{4}dx \lesssim \sup_{[0,T]} \|u(t)\|_{\dot{H}^{\frac{1}{2}}}^2\|u(t)\|_{L^{2}}^2$$
which is the 3d interaction Morawetz estimate that appears in \cite{ckstt4}.
\begin{rem}
Note that although we start with different solutions $u_{1}$ and $u_{2}$ at the end we specialize to $u_{1}=u_{2}=u$. We will omit this step because the notation can be confused with the abbreviation $u_i=u(x_i)$. The meaning will always be clear
 from the context and thus we avoid to introduce a notation that will read $u_{i}(x_{i}):=u_{i}^{i}$ for different solutions
 taking values in $\Bbb R^{n_{i}}$.
\end{rem} 
\subsection{Interaction Morawetz inequality in two dimensions.}
For $n=2$ (in that case $(x_1,x_2)\in \mathbb R^2 \times \mathbb R^2$) we proceed as follows:
\\
Let $f:[0,\infty) \rightarrow [0,\infty)$ be such that
\[f(x):= \left\{\begin{array}{ll}
\frac{1}{2M}x^{2}(1-\log{\frac{x}{M}})& \mbox{if $|x|<\frac{M}{\sqrt{e}}$}\\
100x & \mbox{if $|x|>M$}\\
smooth\ and\ convex\ for\ all\ x & \mbox{}
\end{array}
\right.\]
and $M$ is a large parameter that we will choose later. It is obvious that the functions $\frac{1}{2M}x^{2}(1-\log{\frac{x}{M}})$
 and $100x$ are convex in their domain, and the graph of either function lies strictly above the tangent lines of the other.
 Thus one can construct a function with the above properties. If we apply Proposition \ref{wmor} with the weight 
$a(x_{1},x_{2})=f(|x_1-x_2|)$ and tensoring again two functions we conclude that
$$\int_0^T\int_{\mathbb R^2 \times \mathbb R^2}(-\Delta \Delta a(x_1,x_2))|u(x_{1},t)|^2|u(x_2,t)|^2dx_1dx_2dt \lesssim 2
\sup_{[0,T]}|M_{a}^{\otimes_{2}}(t)|$$
But for $|x_1-x_2|<\frac{M}{\sqrt{e}}$ we have that $\Delta a(x_1,x_2)= \frac{2}{M}\log (\frac{M}{|x_1-x_2|})$ and thus 
$$-\Delta \Delta a(x_1,x_2)=\frac{2}{M}\delta_{\{x_{1}=x_{2}\}}.$$ On the other hand for $|x_1-x_2|>M$ we have that
 $$-\Delta \Delta a(x_1,x_2)= O(\frac{1}{|x_1-x_2|^3})=O(\frac{1}{M^{3}}).$$
We have a similar bound in the region in between just because
$a(x_1,x_2)$ is smooth, so all in all, we have
$$-\Delta \Delta a(x_1,x_2)=\frac{2}{M}\delta_{\{x_{1}=x_{2}\}}+O(\frac{1}{M^{3}}).$$
Thus 
$$\int_0^T\int_{\mathbb R^2 \times \mathbb R^2}(-\Delta \Delta a(x_1,x_2))|u(x_{1},t)|^2|u(x_2,t)|^2dx_1dx_2dt=
\frac{2}{M}\int_0^T\int_{\mathbb R^2}|u(x,t)|^4dxdt+$$
$$O(\frac{1}{M^{3}})\int_0^T\int_{\mathbb R^2 \times \mathbb R^2}
|u(x_1,t)|^2|u(x_2,t)|^2dx_1dx_2dt.$$
By Fubini's Theorem
\begin{equation}
\label{crudestep}
\frac{C}{M^{3}}\int_0^T\int_{\mathbb R^2 \times \mathbb R^2}
|u(x_1,t)|^2|u(x_2,t)|^2dx_1dx_2dt \lesssim
\frac{CT}{M^{3}}\|u\|_{L_{t}^{\infty}L_{x}^{2}}^4.
\end{equation}
On the other hand
$$\sup_{[0,T]}|M_{a}^{\otimes_{2}}(t)|\lesssim \sup_{[0,T]}\|u\|_{L_{t}^{\infty}L_{x}^{2}}^3\|u\|_{L_{t}^{\infty}\dot{H}_{x}^{1}}.$$
Thus by applying Proposition \ref{wmor}
$$\frac{2}{M}\int_0^T\int_{\mathbb R^2}|u(x,t)|^4dxdt\lesssim 
\sup_{[0,T]}\|u\|_{L_{t}^{\infty}L_{x}^{2}}^3\|u\|_{L_{t}^{\infty}\dot{H}_{x}^{1}}+\frac{CT}{M^{3}}\|u\|_{L_{t}^{\infty}L_{x}^{2}}^4.$$
Multiplying the above equation by $M$ and balancing the two terms on the right hand side 
by picking $M \sim T^{\frac{1}{3}}$ we get a better estimate than was obtained in \cite{fg}
$$\|u\|_{L_{t\in [0,T]}^{4}L_{x}^{4}}^4 \lesssim T^{\frac{1}{3}}
\sup_{[0,T]}\|u\|_{L_{t}^{\infty}L_{x}^{2}}^3\|u\|_{L_{t}^{\infty}\dot{H}_{x}^{1}}+T^{\frac{1}{3}}\|u\|_{L_{t}^{\infty}L_{x}^{2}}^4.$$
We note that $\nabla a \in L^{\infty}(\mathbb R^2)$ an observation that will be used strongly later.
\subsection{A new a priori Strichartz estimate in one dimension}
We can combine the calculations that we did in the two dimensional case with the work in \cite{chvz} and obtain the following 
estimate in one dimension,
\begin{equation}\label{66}
\|u\|_{L_{t\in [0,T]}^{6}L_{x}^{6}}^6 \lesssim T^{\frac{1}{3}}
\sup_{[0,T]}\|u\|_{L_{t}^{\infty}L_{x}^{2}}^5\|u\|_{L_{t}^{\infty}\dot{H}_{x}^{1}}+T^{\frac{1}{3}}\|u\|_{L_{t}^{\infty}L_{x}^{2}}^6.
\end{equation}
We can derive this estimate by considering $(x_1,x_2,x_3) \in \Bbb R \times \Bbb R \times \Bbb R$, tensoring three solutions, and
 using an orthonormal change of variables $z=Ax$ where $A$ is an orthonormal matrix. Then we pick the convex weight function to be
\[a(z_1,z_2,z_3)= \left\{\begin{array}{ll}
\frac{1}{2M}(z_2^2+z_3^2)(1-\log{\frac{(z_2^2+z_3^2)^{1/2}}{M}})& \mbox{if $(z_2^2+z_3^2)^{1/2}<\frac{M}{\sqrt{e}}$}\\
100(z_2^2+z_3^2)^{1/2} & \mbox{if $(z_2^2+z_3^2)^{1/2}>M$}\\
smooth\ and\ convex\ for\ all\ z \in \Bbb R^3 & \mbox{}
\end{array}
\right.\]
But then $\Delta_{z}=\Delta_{x}$ and an explicit calculation shows that
$$\Delta_{z}a(z_1,z_2,z_3)=C\frac{1}{M}\delta_{z_2 = z_3}+O(\frac{1}{M^3}).$$
Thus balancing the two terms as in the two dimensional case and going back to the original variables, for details see \cite{chvz},
we obtain \eqref{66}.
\\
\\
\begin{rem} A similar estimate can be obtained if one interpolates the one dimensional estimate
$$\|u\|_{L_{t\in [0,T]}^{8}L_{x}^{8}}^8 \lesssim \sup_{[0,T]}\|u\|_{L_{t}^{\infty}L_{x}^{2}}^7\|u\|_{L_{t}^{\infty}\dot{H}_{x}^{1}}$$
that was proved in \cite{chvz} and the trivial estimate
$$\|u\|_{L_{t\in [0,T]}^{2}L_{x}^{2}} \lesssim T^{\frac{1}{2}}\|u\|_{L_{t\in [0,T]}^{\infty}L_{x}^{2}} \lesssim 
T^{\frac{1}{2}}\|u_{0}\|_{L_{x}^{2}}$$
where we used H\"older's inequality in time and the conservation of mass. This kind of estimate has been already used
in \cite{tz2} to improve the known global well-posedness results for the quintic defocusing nonlinear Schr\"odinger equation in
 one dimension.
\end{rem}
\subsection{Interaction Morawetz inequality in two dimensions for the $I$-NLS equation and the proof of Theorem \ref{main}.}
We now proceed to prove Theorem \ref{main}. 
For motivational purposes let's consider the solution $Iu$ of 
\begin{equation}\label{Iu}
iIu_{t}+\Delta Iu=I(|u|^2u), ~ (x,t) \in \R^2 \times [0,T].
\end{equation}
If $Iu$ would solve not \eqref{Iu} but the nonlinear Schr\"odinger equation
\begin{equation}
iIu_{t}+\Delta Iu=|Iu|^2Iu, ~ (x,t) \in \R^2 \times [0,T],
\end{equation}
then the calculations that we did above in two dimensions would reveal that
$$\|Iu\|_{L_{t\in [0,T]}^{4}L_{x}^{4}}^4 \lesssim T^{\frac{1}{3}}
\sup_{[0,T]}\|Iu\|_{L_{t}^{\infty}L_{x}^{2}}^3\|Iu\|_{L_{t}^{\infty}\dot{H}_{x}^{1}}+
T^{\frac{1}{3}}\|Iu\|_{L_{t}^{\infty}L_{x}^{2}}^4.$$
Of course this is not the case. But we can rewrite equation \eqref{Iu}
as 
\begin{equation}
iIu_{t}+\Delta Iu=|Iu|^2Iu+I(|u|^2u)-|Iu|^2Iu=F(Iu)+(IF(u)-F(Iu)).
\end{equation}
Then if we repeat the calculations, the commutator $IF(u)-F(Iu)$ will indroduce an error term while the term $F(Iu)$ again gives
 rise to a defocusing potential. Thus by Corollary \ref{cor} we get
$$\int_0^T\int_{\Bbb R^2}(-\Delta \Delta a)|Iu(x,t)|^2dxdt \lesssim 2 \sup_{0,T}|\int_{\Bbb R^2}\nabla a(x)
\cdot \Im(\overline{Iu}(x)\nabla Iu)dx|+$$
$$\left| \int_0^T\int_{\Bbb R^2}\nabla a \cdot
\{IF(u)-F(Iu),Iu(x,t)\}_{p}dxdt\right|.$$
The second term on the right hand side of the inequality is what we
call an {\it{Error}} in \eqref{IuMorawetzThird}. We now turn to the
details. The conjugates will play no crucial role in the upcoming argument.
Let us set
$$I U(x,t) = I\otimes I\left( u(x_{1},t)\otimes u(x_{2},t)\right)=\prod_{j=1}^{2}Iu(x_{j},t).$$  
If $u$ solves (\ref{NLS}) for
$n=2,$ then we observe that $IU$ solves (\ref{NLS}) for
$n=4,$ with right hand side $\widetilde{\CN}_I$ given by
$$\widetilde{\CN}_I=\sum_{i=1}^{2}(I(\widetilde{\CN}_{i})\prod_{j=1,j \ne i}^{2}Iu_{j}).$$
Now let us decompose,
$$\widetilde{\CN}_I=\widetilde{\CN}_{good}+\widetilde{\CN}_{bad}=\sum_{i=1}^{2}(\widetilde{\CN}_{i}(Iu) \prod_{j=1,j \ne i}^{2}Iu_{j})+
\sum_{i=1}^{2}\left(I(\widetilde{\CN}_{i})-\widetilde{\CN}_{i}(Iu_{i})\right)\prod_{j=1,j \ne
i}^{2}Iu_{j}$$ The first summand creates a defocusing potential like in the applications before. Thus after integration
 by parts, it creates a positive term that we can ignore. The term we call $\widetilde{\CN}_{bad}$, 
produces the {\it{Error}} term. Repeating the calculations above with $Iu$ instead of $u$ we have the bound:
\begin{align}\label{amor}
\|Iu\|_{L_{T}^{4}L_{x}^{4}}^{4} &\lesssim
T^{\frac{1}{3}}\sup_{0,T}\|Iu\|_{\dot{H}^{1}}\|Iu\|_{L^{2}}^{3}+T^{\frac{1}{3}}\|u_{0}\|_{L_{t}^{\infty}L_{x}^{2}}^4\\\nonumber
&+T^{\frac{1}{3}}\left| \int_0^T\int_{\Bbb R^4}\nabla a \cdot
\{\widetilde{\CN}_{bad},Iu(x_1,t)Iu(x_2,t)\}_{p}dx_1dx_2dt \right|.
\end{align}
Note that we also used the fact that 
$$\|Iu\|_{L^{2}} \lesssim \|u\|_{L^{2}}=\|u_{0}\|_{L^{2}}$$
which follows by the definition of the $I$-operator and conservation of mass.
Note that the third term of \eqref{amor} 
comes from the momentum bracket term in the proof of Proposition
 \ref{wmor}. We also remark that $\nabla a$ is real valued, thus
$$\nabla a \cdot \Re (f \nabla \bar{g}-g \nabla \bar{f})=
\Re \left( \nabla a \cdot (f \nabla \bar{g}-g \nabla \bar{f})\right)$$
and that $\nabla=(\nabla_{x_{1}},\nabla_{x_{2}})$. 
We now wish to compute the dot
product under the integral in \eqref{amor}, that is
$$\Re\left\{\sum_{i=1}^{2}\nabla_{x_{i}}a
\left(\widetilde{\CN}_{bad}(\nabla_{x_{i}}({\overline{Iu_1}}{\overline{Iu_2}})) -
Iu_1 Iu_2\nabla_{x_{i}}{\overline{\widetilde{\CN}}}_{bad} \right)\right\}.$$ 
We start by computing the first summand. Recall that
$$\widetilde{\CN}_{bad}=\sum_{i=1}^{2}\left(I(\widetilde{\CN}_{i})-\widetilde{\CN}_{i}(Iu_{i})\right)\prod_{j=1,j
  \ne i}^{2}Iu_{j}.$$
Using the definition of
$\widetilde{\CN}_{bad},$ and the fact that $\nabla_{x_{1}}$ acts only on $Iu_1$
a direct calculation shows
$$\widetilde{\CN}_{bad}\left(\nabla_{x_1}({\overline{Iu_1}}{\overline{Iu_2}})\right)-Iu_1
Iu_2 (\nabla_{x_1}{\overline{\widetilde{\CN}}}_{bad})=$$
$$\left[(I(\widetilde{\CN}_1)-\widetilde{\CN}(Iu_1))\nabla_{x_1}{\overline{Iu_1}}-\nabla_{x_1}(I(\widetilde{\CN}_1)-
{\overline{\widetilde{\CN}}}(Iu_1))Iu_1\right]|Iu_2|^2.$$
Hence the first summand is given by,
\begin{equation*}\Re\{\nabla_{x_1}a\left[(I(\widetilde{\CN}_1)-
\widetilde{\CN}(Iu_1))\nabla_{x_1}{\overline{Iu_1}}-\nabla_{x_1}(I({\overline{\widetilde{\CN}}}_1)-
{\overline{\widetilde{\CN}}}(Iu_1))Iu_1\right]|Iu_2|^2\}.
\end{equation*}
\\
Analogously one can see that the second summand is
given by:
$$\Re\{\nabla_{x_2}a\left[(I(\widetilde{\CN}_2)-
\widetilde{\CN}(Iu_2))\nabla_{x_2}{\overline{Iu_2}}-\nabla_{x_2}(I({\overline{\widetilde{\CN}}}_2)-
{\overline{\widetilde{\CN}}}(Iu_2))Iu_2\right]|Iu_{1}|^2\}.$$ Thus,
our error term
\begin{align}\mathcal{E}=\int_{0}^{T}\int_{\mathbb{R}^4}\nabla & a \cdot \left(\widetilde{\CN}_{bad}(\nabla(Iu_1Iu_2)) -
\prod_{j=1}^{2}Iu_j\nabla \widetilde{\CN}_{bad}\right)\nonumber\end{align}
reduces to
\begin{align*}\mathcal{E}=\Re\{\int_{0}^{T}\int_{\mathbb{R}^4}\sum_{i=1}^{2}&\{\nabla_{x_i}a\left[(I(\widetilde{\CN}_i)-
\widetilde{\CN}(Iu_i))\nabla_{x_i}{\overline{Iu_i}}-\nabla_{x_i}(I({\overline{\widetilde{\CN}_i}})-
{\overline{\widetilde{\CN}}}(Iu_i) ) Iu_i\right]\\&\times|\prod_{j=1,j\neq
i}^{j=2}Iu_j|^2\}dx_1dx_2dt\}.\end{align*} Hence, by
symmetry,
\begin{equation}\label{total error}|\mathcal{E}| \lesssim |E|,
\end{equation}
where
\begin{align*}E=\int_{0}^{T}\int_{\mathbb{R}^4}&\{\nabla_{x_1}a\left[(I(\widetilde{\CN}_1)-
\widetilde{\CN}(Iu_1))\nabla_{x_1}{\overline{Iu_1}}-\nabla_{x_1}(I({\overline{\widetilde{\CN}}}_1)-
{\overline{\widetilde{\CN}}}(Iu_1))Iu_1\right]\\&\times|Iu_{2}|^2\}dx_1dx_2dt.\end{align*}

We have,
\begin{equation}\label{E}|E| \leq E_1 + E_2
\end{equation} where

$$E_1= \int_{0}^{T}\int_{\mathbb{R}^4}|\nabla_{x_1}a||I(\widetilde{\CN}_1)-
\widetilde{\CN}(Iu_1)||\nabla_{x_1} Iu_1||Iu_2|^2dx_1dx_2dt$$
and

$$E_2= \int_{0}^{T}\int_{\mathbb{R}^4}|\nabla_{x_1}a||\nabla_{x_1}(I(\widetilde{\CN}_1)-
\widetilde{\CN}(Iu_1)||Iu_1||Iu_2|^2dx_1dx_2dt.$$
Since $|\nabla_{x_1}a| \lesssim 1$ applying Fubini's theorem we have

$$E_1 \leq \left(\int_{0}^{T}\int_{\mathbb R^2}|I(\widetilde{\CN}_1)-
\widetilde{\CN}(Iu_1)||\nabla_{x_1}Iu_1|dx_1dt\right)\|Iu\|^2_{L^{\infty}_tL^2_{x}}$$

and

$$E_2 \leq \left(\int_{0}^{T}\int_{\mathbb R^2}|\nabla_{x_1}(I(\widetilde{\CN}_1)-
\widetilde{\CN}(Iu_1))||Iu_1|dx_1dt\right)\|Iu\|^2_{L^{\infty}_tL^2_{x}}.$$

Since the pair $(\infty,2)$ is admissible and by renaming $x_{1}=x$ we have
$$E_1 \leq \left(\int_{0}^{T}\int_{\mathbb R^2}|I(\widetilde{\CN}_1)-
\widetilde{\CN}(Iu_1)||\nabla Iu_1|dxdt\right)Z_{I}^2$$
and
$$E_2 \leq \left(\int_{0}^{T}\int_{\mathbb R^2}|\nabla (I(\widetilde{\CN}_1)-
\widetilde{\CN}(Iu_1))||Iu_1|dxdt\right)Z_{I}^2.$$
Therefore,

$$E_1 \leq \|I(\widetilde{\CN})-
\widetilde{\CN}(Iu)\|_{L^1_tL^2_x}\|\nabla Iu\|_{L^\infty_tL^2_x}Z_I^2$$

and

$$E_2 \leq \|\nabla (I(\widetilde{\CN})-
\widetilde{\CN}(Iu))\|_{L^1_tL^2_x}\|Iu\|_{L^\infty_tL^2_x}Z_I^2.$$

Again, since $(\infty,2)$ is admissible we obtain:

$$E_1 \leq \|I(\widetilde{\CN})-
\widetilde{\CN}(Iu)\|_{L^1_tL^2_x}Z_I^3$$

and

$$E_2 \leq \|\nabla_x(I(\widetilde{\CN})-
\widetilde{\CN}(Iu))\|_{L^1_tL^2_x}Z_I^3.$$

Therefore, from \eqref{E} and the bounds above, we deduce that

\begin{equation}\label{error}|E| \leq \left (\|I(\widetilde{\CN})- \widetilde{\CN}(Iu)\|_{L^1_tL^2_x}+
\|\nabla_x(I(\widetilde{\CN})-\widetilde{\CN}(Iu))\|_{L^1_tL^2_x}\right )Z_I^3.\end{equation}
We proceed to estimate $\|\nabla (I(\widetilde{\CN})- \widetilde{\CN}(Iu))\|_{L^1_tL^2_x}$,
which is the hardest of the two terms. Toward this aim, let us
observe that since $\widetilde{\CN}(u)=|u|^2u$, we will be able to
work on the Fourier side to estimate the commutator $I(\widetilde{\CN})- \widetilde{\CN}(Iu)$.

We compute, \footnote{We ignore complex conjugates, since our
computations are not affected by conjugation}
$$\widehat{\nabla_x(I(\widetilde{\CN})- 
\widetilde{\CN}(Iu))}(\xi) = \int_{\xi=\xi_1+\xi_2+\xi_3} i\xi[m(\xi)-m(\xi_1)m(\xi_2)m(\xi_3)]
\hat{u}(\xi_1)\hat{u}(\xi_2)\hat{u}(\xi_3) d\xi_1d\xi_2d\xi_3$$

We decompose $u$ into a sum of dyadic pieces $u_j$ localized
around $N_j$. Note that the actual decay of the error term is of order $O(\frac{1}{N^{1-\epsilon}})$. This is because we have
 to keep a factor of size $\max_{j=1,2,3}N_{j}^{-\epsilon}$ in order to sum the different Littlewood-Paley pieces. For
 the simplicity of the argument we omit this technicality that doesn't affect the final result. Then,
\begin{align} & \|\nabla(I(\widetilde{\CN})-
\widetilde{\CN}(Iu))\|_{L^1_tL^2_x}=\|\widehat{\nabla(I(\widetilde{\CN})-
\widetilde{\CN}(Iu))}\|_{L^1_tL^2_\xi}\label{fourierside}\\ \nonumber&\leq
\sum_{N_1,N_2, N_3} \|\int_{\xi=\xi_1+\xi_{2}+\xi_3; |\xi_i|\sim
N_i} |\xi|[m(\xi)-m(\xi_1)m(\xi_2)m(\xi_3)] \widehat{u_1}\widehat{u_2}
\widehat{u_3} d\xi_1d\xi_2d\xi_3\|_{L^1_tL^2_\xi}\\
\nonumber&=\sum_{N_1,N_2, N_3} \|\int_{\xi=\xi_1+\xi_2+\xi_3;
|\xi_i|\sim N_i} |\xi|\frac{[m(\xi)-m(\xi_1)m(\xi_2)
m(\xi_3)]}{m(\xi_1)m(\xi_2)m(\xi_3)} \widehat{Iu_1}\widehat{Iu_2}
\widehat{Iu_3} d\xi_1d\xi_2d\xi_3\|_{L^1_tL^2_\xi}.\end{align}

Without loss of generality, we can assume that the $N_j$'s are
rearranged so that $$N_1 \geq N_2 \geq N_3.$$ Set,
$$\sigma(\xi_1,\xi_2, \xi_3)=|\xi_1 +\xi_2+\xi_3|\frac{[m(\xi_1 +\xi_2+\xi_3)-m(\xi_1)m(\xi_2)m(\xi_3)]}
{m(\xi_1)m(\xi_2)m(\xi_3)}.$$ Then,

$$\sigma(\xi_1,\xi_2,\xi_3)= \sum_{j=1}^{4}\chi_j(\xi_1,\xi_2,\xi_3) \sigma(\xi_1,\xi_2,\xi_3)=
\sum_{j=1}^{4}\sigma_j(\xi_1,\xi_2,\xi_3),$$ where
$\chi_j(\xi_1,\xi_2,\xi_3)$ is a smooth characteristic function
of the set $\Omega_j$ defined as follows:
\begin{itemize}
\item $\Omega_1=\{|\xi_i| \sim N_i, i=1,2,3 ; N_1 \ll N\}$.
\item$\Omega_2=\{|\xi_i| \sim N_i, i=1,2,3 ;N_1 \gtrsim N \gg
N_2\}.$ 
\item $\Omega_3=\{|\xi_i| \sim N_i, i=1,2,3 ;N_1 \geq
N_2 \gtrsim N \gg N_3.\}$ 
\item $\Omega_4=\{|\xi_i| \sim N_i,
i=1,2,3 ;N_1 \geq N_2 \geq N_3 \gtrsim N\}.$ 
\end{itemize}

Hence, from \eqref{fourierside} we get,

\begin{align}&\label{bound}\|\nabla(I(\widetilde{\CN})-
\widetilde{\CN}(Iu))\|_{L^1_tL^2_x}\\ &\lesssim \sum_{N_1,N_2, N_3}
\sum_{j=1}^{4}\|\int_{\xi=\xi_1+\xi_2+\xi_3}
\sigma_j(\xi_1,\xi_2,\xi_3)\widehat{Iu_1}\widehat{Iu_2}
\widehat{Iu_3} d\xi_1d\xi_2d\xi_3\|_{L^1_tL^2_\xi}
=\sum_{N_1,N_2,N_3} \sum_{j=1}^{4}L_j\nonumber.\end{align}

We proceed to analyze the contribution of each of the integrals
$L_j.$

\textbf{Contribution of $L_1$.} Since $\sigma_1$ is identically
zero when $N \geq 4 N_1$, $L_1$ gives no contribution to the sum above.

\textbf{Contribution of $L_2$.} We have,

$$\|\int_{\xi=\xi_1+\xi_2+\xi_3}
\sigma_2(\xi_1+\xi_2+\xi_3) \widehat{Iu_1}\widehat{Iu_2}
\widehat{Iu_3} d\xi_1d\xi_2d\xi_3\|_{L^1_tL^2_\xi}
$$
$$=\frac{1}{N}\|\int_{\xi=\xi_1+\xi_2+\xi_3}
\frac{N}{\xi_1\xi_2}\sigma_2(\xi_1,\xi_2,\xi_3)\widehat{\nabla
Iu_1}\widehat{\nabla Iu_2}\widehat{Iu_3} d\xi_1d\xi_2d\xi_3\|_{L^1_tL^2_\xi}$$
$$\lesssim \frac{1}{N} \|\nabla Iu_1\|_{L^3_tL^{6}_x}
\|\nabla Iu_2\|_{L^3_tL^{6}_x}\|Iu_3\|_{L^3_tL^{6}_x}$$
where in the last line we used the Coifman-Meyer multiplier
theorem, \cite{cm}, and H\"older in time. The application of the multiplier
theorem is justified by the fact that the symbol
$$a_2(\xi_1,\xi_2,\xi_3)=
\frac{N}{\xi_1\xi_2}\sigma_2(\xi_1,\xi_2,\xi_3)$$ is of order
zero. The $L^\infty$ bound follows after an application of the
mean value theorem. Indeed,
$$|a_2(\xi_1,\xi_2,\xi_3)| \leq \frac{N}{N_1N_2}|\xi_1+\xi_2+\xi_3|\frac{|\nabla m(\xi_1)(\xi_2+\xi_3)
|}{m(\xi_1)} \lesssim \frac{N}{N_1N_2} N_1 \frac{N_2}{N_1}\lesssim
1.$$

\textbf{Contribution of $L_3$.} We have,

$$\|\int_{\xi=\xi_1+\xi_2+\xi_3}
\sigma_3(\xi_1+\xi_2+\xi_3) \widehat{Iu_1}\widehat{Iu_2}
\widehat{Iu_3} d\xi_1 d\xi_2 d\xi_3\|_{L^1_tL^2_\xi}$$
$$=\frac{1}{N}\|\int_{\xi=\xi_1+\xi_2+\xi_3}
\frac{N}{\xi_1\xi_2}\sigma_3(\xi_1,\xi_2,\xi_3)\widehat{\nabla
Iu_1}\widehat{\nabla Iu_2}\widehat{Iu_3} d\xi_1
d\xi_2d\xi_3\|_{L^1_tL^2_\xi}$$
$$\lesssim \frac{1}{N} \|\nabla Iu_1\|_{L^3_tL^{6}_x}
\|\nabla Iu_2\|_{L^3_tL^{6}_x}\|Iu_3\|_{L^3_tL^{6}_x}$$
where in the last line we used the Coifman-Meyer multiplier
theorem, and Holder in time. The application of the multiplier
theorem is justified by the fact that the symbol
$$a_3(\xi_1,\xi_2,\xi_3)=
\frac{N}{\xi_1\xi_2}\sigma_3(\xi_1,\xi_2,\xi_3)$$ is of order
zero. The $L^\infty$ bound follows from the following chain of
inequalities,

$$|a_3(\xi_1, \xi_2, \xi_3)| \lesssim \frac{N}{N_1N_2}
|\xi_1+\xi_2+\xi_3|(\frac{m(\xi_1+\xi_2+\xi_3)}{m(\xi_1)m(\xi_2)}+1)
$$$$\lesssim  \frac{N}{N_1N_2}(\frac{N_1}{m(N_2)} + N_1) \lesssim 1,
$$ where we have used the fact that $|\xi|m(\xi)$ is monotone increasing for any $s>0$ and thus 
$$|(\xi_1+\xi_2+\xi_3)m(\xi_1+\xi_2+\xi_3)|\lesssim |\xi_{1}m(\xi_{1})|.$$

\textbf{Contribution of $L_4$.}
We continue as above using the $L^3_tL^{6}_x$ Strichartz norms and get
$$\|\int_{\xi=\xi_1+\xi_2+\xi_3}
\sigma_4(\xi_1+\xi_2+\xi_3) \widehat{Iu_1}\widehat{Iu_2}
\widehat{Iu_3} d\xi_1d\xi_2 d\xi_3\|_{L^1_tL^2_\xi}\lesssim
\frac{1}{N^2} \prod_{j=1}^{3}\|\nabla Iu_j\|_{L^3_tL^{6}_x},$$
where in this case the symbol to which we apply the multiplier
theorem is:

$$a_4(\xi_1,\xi_2,\xi_3)=
\frac{N^2}{\xi_1\xi_2\xi_3}\sigma_4(\xi_1,\xi_2,\xi_3).$$

In all the cases above, we proved the $L^\infty$ bound for the
symbols $a_i(\xi_1,\xi_2,\xi_3), i=2, 3, 4.$ The reader can easily verify the conditions of the Coifman-Meyer theorem
 for the higher order derivatives.











Finally, since the pair $(3,6)$ is admissible, we obtain that in
all of the cases above

$$\|\int_{\xi=\xi_1+\xi_2+\xi_3}
\sigma_i(\xi_1+\xi_2+\xi_3) \widehat{Iu_1} \widehat{Iu_2}
\widehat{Iu_3} d\xi_1 d\xi_2d\xi_3\|_{L^1_tL^2_\xi}\lesssim \frac{1}{N}Z_I^3 .$$

Therefore, we deduce from \eqref{bound} that
$$\|\nabla(I(\widetilde{\CN})- \widetilde{\CN}(Iu))\|_{L^1_tL^2_x} \lesssim
\frac{1}{N^{1-}}Z_I^3.$$

Analogously, $$\|I(\widetilde{\CN})- \widetilde{\CN}(Iu)\|_{L^1_tL^2_x} \lesssim
\frac{1}{N^{1-}}Z_I^3.$$ Hence, in view of \eqref{error} we obtain
the following estimate for the error term,

$$|E| \lesssim \frac{1}{N^{1-}}Z_I^{6}.$$
Thus, \eqref{total error} implies
$$\left|\int_{0}^{T}\int_{\mathbb{R}^4}\nabla a \cdot \left(\widetilde{\CN}_{bad}(\nabla(Iu_1Iu_2) -
Iu_1Iu_2\nabla \widetilde{\CN}_{bad}\right)\right| \lesssim
\frac{1}{N^{1-}}Z_I^{6}.$$
This completes the proof of Theorem \ref{main}.

\section{Proof of the main Theorem and comments on further refinements.}
\subsection{Proof of the main Theorem.}
\begin{proof}
Suppose that $u(t,x)$ is a global in time solution to
\eqref{ivp1} with  initial data $u_0 \in C_0^\infty(\mathbb R^2)$. We
will prove that $\|u(t)\|_{H^s}$ obeys polynomial-in-time upper bounds
with the implied constants not depending upon the extra decay and
regularity properties of $u_0$. A familiar density argument then
establishes that $\eqref{ivp1}$ is globally well-posed for $H^s$
initial data in the range of $s$ for which we prove the polynomial
bounds, namely for $s > \frac{2}{5}$. Set
$u^\lambda(x)=\frac{1}{\lambda}u(\frac{x}{\lambda},\frac{t}{\lambda^{2}})$.
We choose the parameter $\lambda$ so that $\|I u_0^\lambda\|_{H^1}
= O(1)$, that is
\begin{equation}\label{L}
\lambda \sim N^{\frac{1-s}{s}}.
\end{equation}

Next, let us pick a time
$T_0$ arbitrarily large, and let us define
\begin{equation}
S : = \{0 < t < \lambda^2T_0 :
\|Iu^\lambda\|_{L^{4}_tL^{4}_x([0,t]\times
\mathbb R^2 )} \leq KN^{\frac{1}{8}}t^{\frac{1}{12}}\},
\end{equation}
with $K$ a constant to be chosen later. Notice that this choice of the set $S$ is dictated by the apriori estimate
of the $L_t^4L_x^4$ norm of $Iu$. Look also equation \eqref{L4error}.

We claim that $S$ is the
whole interval $[0,\lambda^2T_0].$ Indeed, assume by contradiction
that it is not so, then since
$$\|Iu^\lambda\|_{L^{4}_tL^{4}_x([0,t]\times
\mathbb R^2 )}$$ is a continuous function of time, there exists a time $T
\in [0,\lambda^2T_0]$ such that
\begin{align}
\label{contr1}
& \|Iu^\lambda\|_{L^{4}_tL^{4}_x([0,T]\times
\mathbb R^2 )} >KN^{\frac{1}{8}}t^{\frac{1}{12}}\\
\label{contr2}
&\|Iu^\lambda\|_{L^{4}_tL^{4}_x([0,T]\times
\mathbb R^2 )} \leq 2KN^{\frac{1}{8}}t^{\frac{1}{12}}.
\end{align}
We now split the interval $[0,T]$ into subintervals
$J_k$, $k=1,...,L$ in such a way
that
\begin{equation}
\|Iu^\lambda\|^4_{L^{4}_tL^{4}_x([0,J_{k}]\times
\mathbb R^2 )} \leq \mu_0, 
\end{equation} 
with $\mu_0$ as in Proposition \ref{lwp}. This is possible because of \eqref{contr2}. Then, the
number $L$ of possible subintervals must satisfy
\begin{equation} \label{Lc}
L \sim \frac{(2KN^{\frac{1}{8}}T^{\frac{1}{12}})^{4}}{\mu} \sim \frac{(2K)^{4}N^{\frac{1}{2}}T^{\frac{1}{3}}}{\mu}.
\end{equation}
From Proposition \ref{lwp} and Proposition \ref{ac}, we know that,
for any $\frac{1}{3}<s<\frac{1}{2}$
\begin{equation}\label{energybound}
\sup_{[0,T]}E(Iu^\lambda(t)) \lesssim E(Iu_0^\lambda) +
\frac{L}{N^{\frac{3}{2}}}
\end{equation}
and by our choice of $\lambda$, $E(Iu_0^\lambda)\lesssim 1.$
Note that $\frac{2}{5}>\frac{1}{3}$ and we can apply the previous Propositions. Hence, in order to guarantee that
\begin{equation}\label{energyb}
E(Iu^\lambda)\lesssim 1\end{equation} holds for all $t \in [0,T]$
we need to require  that
$$ L \lesssim N^{\frac{3}{2}}. $$
Since $T \leq \lambda^2T_0,$ according to \eqref{L}, this is fulfilled
as long as
\begin{equation} \label{LN}
\frac{(2K)^{4} N^{\frac{1}{2}}(\lambda^2
T_0)^{\frac{1}{3}}}{\mu_0} \sim N^{\frac{3}{2}},
\end{equation}
From our choice \eqref{L} of $\lambda$, the expression \eqref{LN} implies that
\begin{equation}\label{Tn}
T_{0}^{\frac{1}{3}}\frac{(2K)^4}{\mu_0} \sim N^{\frac{5s-2}{3s}}.
\end{equation}
If $s>\frac{2}{5}$, we have that $T_{0}$ is arbitrarily large if we send $N$ to infinity.
\\
\\
We now use the energy control in the Morawetz estimate to show that
\eqref{contr1} fails to hold. Recall the apriori estimate \eqref{amor}
\begin{align}
\|Iu\|_{L_{T}^{4}L_{x}^{4}}^{4} &\lesssim
T^{\frac{1}{3}}\sup_{[0,T]}\|Iu\|_{\dot{H}^{1}}^{1}\|Iu\|_{L^{2}}^{3}+T^{\frac{1}{3}}\|u_{0}\|_{L^{2}}^{4}\\\nonumber
&T^{\frac{1}{3}}\int_0^T\int_{\mathbb R^4}\nabla a \cdot
\{\widetilde{\CN}_{bad},Iu_1Iu_2\}_{p}dx_1dx_2dt.
\end{align}
Let's define
$$Error(t):=\int_{\mathbb R^4}\nabla a \cdot
\{\widetilde{\CN}_{bad},Iu_1Iu_2\}_{p}dx_1dx_2.$$
By Theorem \ref{main} and Proposition \ref{lwp} we know that on each interval $J_{k}$ we have that
$$\int_{J_{k}}Error(t)dt \lesssim \frac{1}{N}Z_{I}^{6} \lesssim \frac{1}{N}\|I
u^\lambda(a)\|_{H^1}^{6} \lesssim \frac{1}{N}$$
and summing all the $J_{k}$'s we have that
$$\int_{0}^{T}Error(t)dt \lesssim L\frac{1}{N} \sim \frac{N^{\frac{3}{2}}}{N} \sim N^{\frac{1}{2}}.$$
But then we have
\begin{equation}\label{L4error}
\|Iu^{\lambda}\|_{L_{T}^{4}L_{x}^{4}}^{4} \lesssim
T^{\frac{1}{3}}\sup_{[0,T]}\|Iu^{\lambda}\|_{\dot{H}^{1}}^2\|Iu^{\lambda}\|_{L^{2}}^{2}+T^{\frac{1}{3}}\|u^{\lambda}\|_{L^{2}}^{4}
\end{equation}
$$+T^{\frac{1}{3}}\int_{0}^{T}Error(t)dt 
\lesssim T^{\frac{1}{3}}+T^{\frac{1}{3}}N^{\frac{1}{2}}\lesssim N^{\frac{1}{2}}T^{\frac{1}{3}}.$$

This estimate contradicts \eqref{contr1} for an appropriate
choice of $K$. Hence $S = [0,\lambda^2T_0]$, and $T_0$ can
be chosen arbitrarily large. In addition, we have also proved that for $s>\frac{2}{5}$
$$\| I u^\lambda(\lambda^{2}T_{0})\|_{H^1_x}=O(1).$$
But then,
$$\|u(T_{0})\|_{H^{s}} \lesssim \|u(T_{0})\|_{L^{2}}+\|u(T_{0})\|_{\dot{H}^{s}}=
\|u_{0}\|_{L^{2}}+\lambda^{s}\|u^{\lambda}(\lambda^{2}T_{0})\|_{\dot{H}^{s}}$$ 
$$\lesssim \lambda^{s}\| I u^\lambda(\lambda^{2}T_{0})\|_{H^1_x}\lesssim \lambda^{s}\lesssim N^{1-s}
\lesssim T_{0}^{\frac{3s(1-s)}{2(5s-2)}}$$
Since $T_{0}$ is arbitrarily
 large, the apriori bound on the $H^{s}$ norm concludes the global well-posedness of the the Cauchy problem \eqref{ivp1}.
\end{proof}


\begin{thebibliography}{100}

\bibitem{jb1}
J. Bourgain, {\it Refinements of Strichartz' inequality and applications to
2D-NLS with critical nonlinearity,}
International Mathematical Research Notices, 5 (1998), 253--283.

\bibitem{jb2}
J. Bourgain.
{\it Global solutions of nonlinear {S}chr\"odinger equations,} American Mathematical Society, Providence, RI, 1999.

\bibitem{tc} T. Cazenave,
{\it Semilinear Schr\"oodinger equations,} CLN 10, eds: AMS, 2003.

\bibitem{cm}  R. Coifman, and Y. Meyer,
{\it Au dela des operateurs pseudo-differentiels,} Asterisque, 57, Societe Mathematique de France, Paris 1978.

\bibitem{chvz} J. Colliander, J. Holmer, M. Visan and X. Zhang,
{\it Global existence and scattering for rough solutions to generalized
nonlinear Schr\"odinger equations on $\mathbb R$,}
Preprint (2006). ({\tt{math.AP/0612452}})

\bibitem{crsw} J. Colliander, S. Raynor, C. Sulem, J. D. Wright, 
{\it Ground state mass concentration in the $L^2$-critical nonlinear
Schr\"odinger equation below $H^1$,} Math Res. Lett. 12 (2005) no. 2-3, 357-375.  


\bibitem{ckstt2} J. Colliander, M. Keel, G. Staffilani, H. Takaoka, T. Tao,
{\it Almost conservation laws and global rough solutions to a nonlinear Schr\"odinger
equation}, Math. Research Letters, 9 (2002), 659-682.


\bibitem{ckstt4} J. Colliander, M. Keel, G. Staffilani, H. Takaoka and
T. Tao,
 {\it Global existence and scattering for rough solutions to a nonlinear Schr\"odinger
 equations on $\R^{3}$ }, C.P.A.M. \textbf{57} (2004), no. 8, 987--1014.


\bibitem{ckstt5} J. Colliander, M. Keel, G. Staffilani, H. Takaoka and T. Tao,
{\it Global well-posedness and scattering in the energy space for the critical nonlinear Schr\"odinger equation in $\mathbb R^{3}$,}
to appear in Annals Math. ({\tt{math.AP/0402129}})

\bibitem{CKSTT:Angles} J. Colliander, M. Keel, G. Staffilani,
  H. Takaoka and T. Tao, {\it {Global well-posednesss for the cubic
    nonlinear Schr\"odinger equation in $H^s (\R^2)$ for $s > 1/2$}},
  preprint, 2007.



\bibitem{tz2} D. De Silva, N. Pavlovic, G. Staffilani, and N. Tzirakis
{\it  Global well-posedness and polynomial bounds for the defocusing
$L^{2}$-critical nonlinear Schr\"odinger equation in $\R$,} submitted
for publication. ({\tt{math.AP/0702707 }})


\bibitem{fg}
Y. Fang and M. Grillakis, {\it On the global existence of rough
solutions of the cubic defocusing Schr\"odinger equation in
$\R^{2+1}$,} to appear in JHDE.

\bibitem{gv} J. Ginibre and G. Velo, 
{\it The global Cauchy problem for the nonlinear Schr\"oodinger equation,} H. Poincar\'e
Analyse Non Lin\'eaire, 2 (1985), 309-327.

\bibitem{kt} M. Keel and T. Tao,
{\it Endpoint Strichartz estimates,} H. Poincar\'e Analyse non
Lin\'eaire, 120 (1998), 955-980.

\bibitem{LinStrauss} J. E. Lin and W. A. Strauss, {\it Decay and
    scattering of solutions of a nonlinear Schr\"odinger equation},
  J. Funct. Anal., \textbf{30}:2, (1978), 245--263.

\bibitem{es} E. M. Stein,
{\it Harmonic Analysis: Real variable Methods, Orthogonality and Oscillatory integrals,} Princeton Univ. Press, Princeton (1993).



\bibitem{tt} T. Tao,
{\it Nonlinear dispersive equations. Local and global analysis} CBMS 106, eds: AMS, 2006.

\bibitem{tvz1}  T. Tao, M. Visan, and X. Zhang
{\it Minimal-mass blowup solutions of the mass-critical NLS,}
preprint, 2006. ({\tt{math.AP/0609690}})

\bibitem{tvz}  T. Tao, M. Visan, and X. Zhang
{\it Global well-posedness and scattering for the mass-critical nonlinear
Schr\"odinger equations for radial data in high dimensions,} preprint,
2006. ({\tt{math.AP/0609692 }})

\end{thebibliography}
\end{document}